\documentclass[12pt,epsf,epic,eepic,float,psfig,fullpage]{article}

\newtheorem{Def}{{\mbox{$\;\;\;\;\;\,$}}Definition}[section]
\newtheorem{Th}{{\mbox{$\;\;\;\;\;\,$}}Theorem}[section]

\newtheorem{Lem}{{\mbox{$\;\;\;\;\;\,$}}Lemma}[section]

\begin{document}

\title{Stability of a Time-varying Fishing Model   \\
with Delay}

\author{
   L. Berezansky
   \thanks{Research supported in part by the Israeli
   Ministry of Absorption.}  \\
   Department of Mathematics and Computer Science  \\
   Ben-Gurion University of Negev  \\
   Beer-Sheva 84105, Israel  \\
   email:  {\tt{brznsky@cs.bgu.ac.il}} \\\\
   L. Idels
   \thanks{Research supported in part by a grant from the Natural
   Sciences and Engineering
   Research Council of Canada.}
      \\
   Department of Mathematics \\
   Malaspina University-College  \\
   900 Fifth St. Nanaimo, BC, Canada \, V9S5J5   \\
   email:  {\tt{idelsl@mala.bc.ca}}  \\\\
   Corresponding Author:  Lev Idels {\vspace{-1mm}}
      }

\maketitle \vfill{}

{\bf Abstract} We introduce a delay differential equation model
which describes how fish are harvested
$$
\dot{N}(t)=\left[\frac{a(t)}{1+\left(\frac{N(\theta
(t))}{K(t)}\right) ^{\gamma }}-b(t)\right]
N(t)~~~~~~~~~~~~~~~~~~~~~~~~~~~~~~~~~~~~~~~~~~~~~~~~(A)
$$
In our previous studies we investigated the persistence of equation
(A) and existence of a periodic solution for this equation.

Here we study the stability (local and global) of the periodic
solutions of equation (A).

\vspace{5mm}

Keywords-Fishery, Periodic Environment, Delay Differential
Equations, Global and Local Stability. \vspace{5mm}

\section{Introduction and Preliminaries}

Consider the following differential equation which is widely used in
Fisheries \cite{Kot,BBI}
\begin{equation}
\label{1} \dot{N}=[ \beta (t,N)-M(t,N)] N,
\end{equation}
where $N=N(t)$ is the population biomass, $\beta (t,N)$ is the
per-capita fecundity rate, and  $M(t,N)$  is the per-capita fishing
mortality rate due to natural mortality causes and harvesting.

In equation (\ref{1}) let $\beta (t,N)$ be a Hill's type function
\cite{Kot, BBI,BI}\
\begin{equation}
\label{2} \beta (t,N)=\frac{a}{1+\left(\frac{N}{K}\right)^{\gamma
}},
\end{equation}
where $a$ and $K$ are positive constants, and  $\gamma >0$ is a
parameter.

Traditional Population Ecology is based on the concept that carrying
capacity does not change over time even though it is known
\cite{Myers} that the values of carrying capacity related to the
habitat areas might vary, e.g., year-to-year changes in weather
affect fish population.

We assume that in (\ref{2}) $a=a(t)$ , $K=K(t),$ and $M(t,N)=b(t)$
are continuous positive functions.

Generally, Fishery models \cite{Kot,BBI} recognize that for real
organisms it takes time to develop from newborns to reproductively
active adults.

Let in equation (\ref{2}) $N=N(\theta (t)),$ where $\theta (t)$ is
the maturation time delay $0\leq \theta (t)\leq t.$ If we take into
account that delay, then we have the following time-lag model based
on equation (\ref{1})
\begin{equation}
\label{3} \dot{N}(t)=\left[\frac{a(t)}{1+\left( \frac{N(\theta
(t))}{K(t)}\right) ^{\gamma }}-b(t)\right] N(t)
\end{equation}
for $\gamma >0$, with the initial function and the initial value
\begin{equation}
\label{4} N(t)=\varphi (t),~t<0,~N(0)=N_{0}
\end{equation}
under the following conditions:

(a1) $a(t),b(t),K(t)$ are continuous on $[ 0,\infty )$ functions,
$b(t)\geq b>0,~K\geq K(t)\geq k>0$;

(a2) $\theta (t)$ is a continuous function, $\theta (t)\leq
t,~\limsup\limits_{t\rightarrow \infty }\theta (t)=\infty $;

(a3) $\varphi :(-\infty ,0)\rightarrow R$ is a continuous bounded
function, $ \varphi (t)\geq 0,N_{0}>0$.

\begin{Def}  \label{globsoln}
{\rm{ A function $N:R\rightarrow R$ with continuous derivative is
called {\em a (global) solution} of problem (\ref{3}), (\ref{4}), if
it satisfies equation (\ref{3}) for all $t\in [0,\infty )$ and
equalities (\ref{4}) for $t\leq 0$. }}
\end{Def}
\noindent If $t_{0}$ is the first point, where the solution $N(t)$
of (\ref{3}), (\ref{4}) vanishes, i.e., $N(t_{0})=0$, then we
consider the only positive solutions of the problem (\ref{3}),
(\ref{4}) on the interval $\lbrack 0,t_{0})$.

Recently \cite{BI} we proved the following results:
\begin{Lem}  \label{Lemma1}
Suppose $a(t)>b(t),$
$$
\sup_{t>0}\int_{\theta (t)}^t (a(s)-b(s))ds <\infty,~
 \sup_{t>0}\int_{\theta (t)}^t b(s)ds <\infty.
$$
Then there exists the global positive solution of (\ref{3}),
(\ref{4}) and this solution is persistence:
$$ 0<\alpha_N\leq N(t)\leq\beta_N<\infty.$$
\end{Lem}

\begin{Lem}   \label{Lemma2}
Let $a(t),b(t),K(t),\theta (t)$  be T-periodic functions, $a(t)\geq
b(t)$. If at least one of the following conditions hold:

(b1) $$\inf\limits_{t\geq
0}\left(\frac{a(t)}{b(t)}-1\right)K^{\gamma }(t)>1,
$$(b2) $$ \sup\limits_{t\geq
0}\left(\frac{a(t)}{b(t)}-1\right)K^{\gamma }(t)<1,
$$
then equation (\ref{3}) has at least one periodic positive solution
$N_{0}(t).$
\end{Lem}

In what follows, we use a classical result from the theory of
differential equations with delay $\cite{Kuang,Gopal}$.

\begin{Lem}  \label{Lemma3}
Suppose that for linear delay differential equation
\begin{equation}
\label{5} \dot{x}(t)+r(t)x(h(t))=0
\end{equation}
where $0\leq t-h(t)\leq\sigma$, the following conditions hold:

\begin{equation}
\label{6}
r(t)\geq r_0>0 ,\\
\end{equation}
\begin{equation}
\label{7} \limsup_{t\rightarrow\infty}\int_{h(t)}^t
r(s)ds<\frac{3}{2}
\end{equation}

Then for every solution $x$ of equation(\ref{5}) we have
$\lim_{t\rightarrow\infty} x(t)=0$.
\end{Lem}

\section{Main Results}

Let us study global stability of the periodic solutions of equation
(\ref{3}).

\begin{Th}  \label{Theorem1}
Let $a(t), b(t),K(t),\theta (t)$  be T-periodic functions,
satisfying conditions of Lemma \ref{Lemma1} and one of  conditions
b1) or b2) of Lemma \ref{Lemma2}. Suppose also that
\begin{equation}
\label{8} a(t)\geq a_0>0,~~\gamma \int\limits_{\theta
(t)}^{t}a(s)ds<6.
\end{equation}
Then there exists the unique positive periodic solution $N_0(t)$ of
(3) and for every positive solution N(t) of (3) we have
$$lim_{t\rightarrow\infty}(N(t)-N_{0}(t)) =0,$$ i.e., the positive
periodic solution $N_{0}(t)$ is a global attractor for all positive
solutions of (\ref{3}).
\end{Th}
{\em{Proof.}} Lemma \ref{Lemma2} implies that there exists a
positive periodic solution $N_0(t)$. If that solution is an
attractor for all positive solutions then it is the unique positive
periodic solution.

We set $N(t)=\exp (x(t))$ and rewrite equation (3) in the form
\begin{equation}
\label{9} \dot{x}(t)=\frac{a(t)}{1+\left( \frac{e^{x(\theta
(t))}}{K(t)} \right) ^{\gamma }}-b(t).
\end{equation}
Suppose $u(t)$ and $v(t)$ are two different solutions of  (\ref{9}).
Denote $w(t)=u(t)-v(t)$. To prove the Theorem \ref{Theorem1} it is
sufficient to show that $\lim_{t\rightarrow\infty} w(t)=0$.

It follows
\begin{equation}
\label{10} \dot{w}(t)=a(t)\left[ \frac{1}{1+\left( \frac{e^{u(\theta
(t))}}{ K(t)}\right) ^{\gamma }}-\frac{1}{1+\left( \frac{e^{v(\theta
(t))}}{K(t)}\right) ^{\gamma }}\right)
\end{equation}
Let
\begin{equation}
\label{11} f(y,t)=\frac{1}{1+\left( \frac{e^{y}}{K(t)}\right)
^{\gamma }}.
\end{equation}
Using the mean value theorem, we have for every $t$
\begin{equation}
\label{12} f(y,t)-f(z,t)=f^{\prime} (c)(y-z),
\end{equation}
where $\min\{y,z\}\leq c(t)\leq \max\{y,z\}$.

Clearly,
\begin{equation}
\label{13} f^{\prime }_y(y,t)=-\frac{\gamma \left(
\frac{e^{y}}{K(t)}\right)  ^{\gamma }}{ \{1+\left(
\frac{e^{y}}{K(t)}\right)  ^{\gamma }\}^{2}}
\end{equation}
and $\left| f^{\prime }_y(y,t)\right| <\frac{1}{4}\gamma .$

Equalities  (\ref{11})- (\ref{12}) imply that equation  (\ref{10})
takes the form
\begin{equation}
\label{14} \dot{w}(t)=-M(t)w(\theta (t)),
\end{equation}
where
$$
M(t)=\frac{\gamma a(t)\left( \frac{e^{c(t)}}{K(t)}\right) ^{\gamma
}}{ \{1+\left( \frac{e^{c(t)}}{K(t)}\right) ^{\gamma }\}^{2}},
$$
and$$ \min\{u(\theta (t)),v(\theta (t))\}\leq c(t)\leq
\max\{u(\theta (t)),v(\theta (t))\}.$$

Now we want to check that for equation(\ref{14}) all conditions of
Lemma \ref{Lemma3} hold.

From (\ref{13})  we have $M(t)<\frac{1}{4}\gamma a(t).$ Therefore
inequality  (7) holds. Let us check inequality (6). Set
$N_1(t)=e^{u(t)},~ N_2(t)=e^{v(t)}$, where $N_1(t), N_2(t)$ are two
solutions of equation (\ref{3}), corresponding to the solutions
$u(t)$ and $v(t)$ of equation (\ref{9}). Lemma \ref{Lemma1} implies
that
$$
M(t)\geq \frac{\gamma a_0
\left(\frac{\min\{\alpha_{N_1},\alpha_{N_2}\}}{K}\right)^{\gamma
}}{\left(1+\left(\frac{\max\{\beta_{N_1},\beta_{N_2}\}}{k}\right)
^{\gamma}\right)^2}>0,
$$
where $\alpha_{N}$ and  $ \beta_{N}$ are defined by Lemma
\ref{Lemma1}. Hence inequality (6) holds and therefore Theorem
\ref{Theorem1} is proven.

Consider now  equation (\ref{3}) with proportional coefficients:
\begin{equation}
\label{15} \dot{N}(t)=\left[ \frac{ar(t)}{1+\left( \frac{N(\theta(t)
)}{K}\right) ^{\gamma }}-br(t)\right) N(t),
\end{equation}
where $r(t)\geq r_0>0$. Clearly, if $a>b$ then equation (\ref{15})
has the unique positive equilibrium
\begin{equation}
\label{16} N^{\ast }=\left( \frac{a}{b}-1\right) ^{\frac{1}{\gamma
}}K.
\end{equation}
{\bf Corollary 1}. If $a>b, \left( \frac{a}{b}-1\right) K^{\gamma
}\neq 1,r(t)\geq r_0>0,$ and
\begin{equation}
\label{17} \gamma a\limsup_{t\rightarrow\infty}\int_{\theta (t)}^t
r(s)ds <6,
\end{equation}
then the equilibrium $N^{\ast }$ is a global attractor for all
positive solutions of equation (\ref{15}).

Let us now compare the global attractivity condition (\ref{17}) with
the local stability conditions.
\begin{Th}\label{Theorem2}
Suppose $a>b, r(t)\geq r_0>0 $ and
\begin{equation}
\label{18} \frac{\gamma
(a-b)b}{a}\limsup_{t\rightarrow\infty}\int_{\theta (t)}^t r(s)ds
<\frac{3}{2}.
\end{equation}
Then the equilibrium $N^{\ast }$ of equation(\ref{15}) is locally
asymptotically stable. \end{Th} {\bf Proof.} Set $x=N-N^{\ast }$ and
from equation (\ref{15}) we have
\begin{equation}
\label{19} \dot{x}(t)=\left[ \frac{ar(t)}{1+\left( \frac{x(\theta(t)
)+N^{\ast }}{K} \right) ^{\gamma }}-br(t)\right] (x(t)+N^{\ast }).
\end{equation}
Denote
$$
F(u,v)=\left[ \frac{ar(t)}{1+\left(  \frac{u+N^{\ast }}{K}\right) ^{\gamma }}%
-br(t)\right] (v+N^{\ast }).
$$
Clearly,
$$
F_{u}^{\prime }(0,0)=-\frac{\gamma (a-b)b}{a}r(t)
$$
and $F_{v}^{\prime }(0,0)=0 $. Hence for equation (\ref{15}) the
linearized equation has a form
\begin{equation} \label{20}
\dot{x}(t)=-\frac{\gamma (a-b)b}{a}r(t)x(\theta(t) ).
\end{equation}
Lemma {\ref{Lemma3} and condition (\ref{18}) imply  that equation
(\ref{20}) is asymptotically stable, therefore the positive
equilibrium $N^{\ast}$ of equation (\ref{15}) is locally
asymptotically stable.

Compare now Theorem \ref{Theorem1} and Theorem \ref{Theorem2}. We
have $\max\{b(a-b)\}=a/4$. Therefore, if
$$
a\gamma\limsup_{t\rightarrow\infty}\int_{\theta (t)}^t r(s)ds<6,
$$
then equation (\ref{15}) has locally asymptotically stable
equilibrium $N^{\ast}$.

The last condition  does not depend on $b$, and is identical to
condition (\ref{17}) that guarantees the existence of a global
attractor. Therefore in Theorem \ref{Theorem2} we obtained the best
possible conditions for global attractivity for equation (\ref{3}).


\begin{thebibliography}{99}
\bibitem{Kot}
M. Kot, Elements of Mathematical Ecology, Cambr.Univ. Press, 2001.
\bibitem{BBI}
L. Berezansky, E. Braverman, and L. Idels, On Delay Differential
Equations With Hill's Type Growth Rate and Linear Harvesting,
''Computers and Mathematics with Applications'',V.49 (2005), 549-563
\bibitem{Myers}
R. Myers, B. MacKenzie, and K. Bowen, What is the carrying capacity
for fish in the ocean? A meta-analysis of population dynamics of
North Atlantic cod, Can. J. Fish. Aquat. Sci. {58} 1464-1476 (2001).
\bibitem{Kuang}
Y. Kuang, Delay Differential Equations With Applications in
Population Dynamics, Academic Press, Inc (1993).
\bibitem{BI} L. Berezansky and L. Idels, Population Models
With Delay in Dynamic Environment, submitted to the Journal
''Nonlinear Analisys-Real World Applications'', 2005(available at
http:/arxiv.org/ftp/math.DS/papers/0601103).
\bibitem{Gopal}
K. Gopalsamy, Stability and Oscillations in Delay Differential
Equations of Population Dynamics, Kluwer Acad. Publ 1992
\end{thebibliography}
\end{document}